\documentclass[reqno,twoside,12pt]{amsart}

\hoffset - 1.5cm

\usepackage{amsmath,amssymb}
\usepackage{wasysym}
\usepackage{xcolor}

\newtheorem{df}{Definition}[section]
\newtheorem{thm}[df]{Theorem}
\newtheorem{lem}[df]{Lemma}

\newtheorem{rem}[df]{Remark}
\newtheorem{exa}[df]{Example}
\title{Green function to the biharmonic equation on $n$-dimensional sphere}
\author{Ilona Iglewska-Nowak}\address{West Pomeranian University of Technology, Department of Mathematics, al. Piastów 17, PL-70-310 Szczecin, Poland, ORCID 0000--0002--1938--8055 } \email{iiglewskanowak@zut.edu.pl}

\begin{document}

\begin{abstract}
Homogeneous and inhomogeneous biharmonic equation are considered on the $n$-dimensional unit sphere. The Green function is given as a series of Gegenbauer polynomials. In the paper, explicit representations of the Green function are found for some sets of parameters.
\end{abstract}


\date{\today}
\maketitle

\section{Introduction}
Biharmonic equation $\Delta^2u=f$ is applied in elasticity theory fluid mechanics, electromagnetism, and continuum mechanics. In the recent years several researchers have been working on finding its solution on the unit sphere. In~\cite{MG17} homogeneous biharmonic equation with polynomial boundary conditions on the $2$-dimensional unit sphere is solved, and the algorithm requires differentiation of the boundary function, but no integration. Paper~\cite{vlV24} presents a solution to the inhomogeneous biharmonic equation on the $2$-dimensional unit sphere being an approximation by a sequence of the same equation, but with special right-hand-sides, which are shifts  of the Dirac delta function. In~\cite{sG18} the author constructs biharmonic functions on $n$-dimensional unit spheres (as well as hyperbolic spaces).

The aim of the present paper is find a Green function to the inhomogeneous and homogeneous biharmonic equation on the $n$-dimensional unit sphere. The series representation of the Green function has been already known. The main contribution of the present paper is an integral representation of the Green function which allows to find an explicit formula for some sets of parameters.

The paper is organized as follows. In Section~\ref{sec: preliminaries} some preliminary information is given. Section~\ref{sec: biharmonic equation} introduces biharmonic equation and the basic representation of its Green function~$G_a$. In Section~\ref{sec: integral representation} integral representations of~$G_a$ are derived which allow computation of explicit formulae for~$G_a$ in Section~\ref{sec: explicit representations}.

\section{Preliminaries}\label{sec: preliminaries}

A square integrable function~$f$ over the $n$-dimensional unit sphere~$\mathcal S^n\subseteq\mathbb R^{n+1}$, $n\geq 2$, with the rotation-invariant measure~$d\sigma$ normalized such that
$$
\Sigma_n=\int_{\mathcal S^n}d\sigma(x)=\frac{2\pi^{(n+1)/2}}{\Gamma\left((n+1)/2\right)},
$$
can be represented as a Fourier series in terms of the hyperspherical harmonics,
\begin{equation}\label{eq: Fs}
f=\sum_{l=0}^\infty\sum_{k\in\mathcal M_{n-1}(l)}a_l^k(f)\,Y_l^k,
\end{equation}
where $\mathcal M_{n-1}(l)$ denotes the set of sequences $k=(k_0,k_1,\ldots,k_{n-1})$ in $\mathbb N_0^{n-1}\times\mathbb Z$ such that $l\geq k_0\geq k_1\geq\ldots\geq|k_{n-1}|$ and $a_l^k(f)$ are the Fourier coefficients of~$f$. The hyperspherical harmonics of degree~$l$ and order~$k$ are given by
\begin{equation}\label{eq: Ylk}
Y_l^k(x)=A_l^k\prod_{\tau=1}^{n-1}C_{k_{\tau-1}-k_\tau}^{\frac{n-\tau}{2}+k_\tau}(\cos\vartheta_\tau)\sin^{k_\tau}\!\vartheta_\tau\cdot e^{\pm ik_{n-1}\varphi}
\end{equation}
for some constants $A_l^k$. Here, $(\vartheta_1,\dots,\vartheta_{n-1},\varphi)$ are the hyperspherical coordinates of \mbox{$x\in\mathcal S^n$},
\begin{align*}
x_1&=\cos\vartheta_1,\\
x_2&=\sin\vartheta_1\cos\vartheta_2,\\
x_3&=\sin\vartheta_1\sin\vartheta_2\cos\vartheta_3,\\
\dots\\
x_{n-1}&=\sin\vartheta_1\dots\sin\vartheta_{n-2}\cos\vartheta_{n-1},\\
x_n&=\sin\vartheta_1\dots\sin\vartheta_{n-2}\sin\vartheta_{n-1}\cos\varphi,\\
x_{n+1}&=\sin\vartheta_1\dots\sin\vartheta_{n-2}\sin\vartheta_{n-1}\sin\varphi,
\end{align*}
and $\mathcal C_\kappa^K$ are the Gegenbauer polynomials of degree~$\kappa$ and order~$K$.
The set of the hyperspherical harmonics of degree~$l$ is denoted by~$\mathcal H_l$.

The Laplace-Beltrami operator on the sphere $\Delta_{\mathcal S^n}$ is defined by 
\begin{equation*}
\Delta_{\mathcal S^n}f=\sum\limits_{k=1}^{n-1}\left(\prod\limits_{j=1}^k\sin\vartheta_j\right)^{-2}
 (\sin \vartheta_k)^{k+2-n}\frac{\partial}{\partial\vartheta_k}\left(\sin^{n-k}\vartheta_k\frac{\partial f}{\partial\vartheta_k}\right)
 +\left(\prod\limits_{j=1}^k\sin\vartheta_j\right)^{-2}\frac{\partial^2f}{\partial\varphi^2}.
\end{equation*}

It satisfies the following relation:
\begin{equation}\label{eq: recursion LB theta}
\Delta_{\mathcal S^n}f=\frac{1}{\sin^{n-1}\vartheta}\frac{\partial}{\partial\vartheta}\left[\sin^{n-1}\vartheta\cdot\frac{\partial f}{\partial\vartheta}\right]
   +\frac{1}{\sin^2\vartheta}\Delta_{\mathcal S^{n-1}}f,\qquad\vartheta=\vartheta_1.
\end{equation}
Setting $t=x_1=\cos\vartheta$ one can write it also as~\cite[formula~(3.8)]{AtkinsonHan}
\begin{equation}\label{eq: recursion LB t}
\Delta_{\mathcal S^n}f=\frac{1}{(1-t^2)^\frac{n-2}{2}}\frac{\partial}{\partial t}\left[(1-t^2)^\frac{n}{2}\cdot\frac{\partial f}{\partial t}\right]
   +\frac{1}{1-t^2}\Delta_{\mathcal S^{n-1}}f.
\end{equation}
In both cases, $\mathcal S^{n-1}$ denotes the $(n-1)$-dimensional unit sphere with hyperspherical coordinates $(\vartheta_2,\dots,\vartheta_{n-1},\varphi)$.

It is known that the hyperspherical harmonics are the eigenfunctions of $\Delta_{\mathcal S^n}$, i.e.,
\begin{equation*}
\Delta_{\mathcal S^n} Y_l^k=-l(l+n-1)Y_l^k,
\end{equation*}
see \cite[Chapter~II, Theorem 4.1]{Shimakura}.
The relation of $\Delta_{\mathcal S^n}$ and the Laplace operator $\Delta$ is given by
\begin{equation*}
\Delta f = R^{-n}\frac{\partial}{\partial R}\left(R^{n}\frac{\partial f}{\partial R}\right) + \frac{1}{R^2}\Delta_{\mathcal S^n} f,
\end{equation*}
where
$R\geq 0$ is the radius of $x\in \mathbb{R}^n$ in the hyperspherical coordinates, see \cite[Chapter~II, Proposition 3.3]{Shimakura}.

The Laplace operator is commutative with $SO(n+1)$-rotations $\Upsilon$,
\begin{equation*}
\Delta\left[f(\Upsilon x)\right]=(\Delta f)(\Upsilon x),
\end{equation*}
see \cite[Chapter~IX, \S~2, Subsec.~4]{Vilenkin}. Consequently, it follows from \eqref{eq: laplacians} that the same holds for the Laplace-Beltrami operator, see also \cite[Chapter~II, formula (3.15)]{Shimakura}.

Since $\mathcal{S}^n$ is a manifold without boundary, the Green's second surface identity implies that for $f,g$ of class $\mathcal C^2$ the following holds:
\begin{equation}\label{eq: green}
\int_{\mathcal{S}^n} \Delta_{\mathcal S^n} f(x) \cdot g(x)\, d\sigma(x)=\int_{\mathcal{S}^n}  f(x) \cdot \Delta_{\mathcal S^n}g(x)\, d\sigma(x).
\end{equation}

The scalar product in $\mathcal L^2(\mathcal S^n)$ is antilinear in the first variable,
$$
\left<f,g\right>=\frac{1}{\Sigma_n}\int_{\mathcal{S}^n}\overline{f(x)}\,g(x)\,d\sigma(x).
$$
Since~$\Delta_{\mathcal S^n}$ is a linear operator, one has
$$
\overline{\Delta_{\mathcal S^n} f}=\Delta_{\mathcal S^n}\overline{f}
$$
and~\eqref{eq: green} can be also written as
\begin{equation*}
\left<\Delta_{\mathcal S^n}f,g\right>=\left<f,\Delta_{\mathcal S^n}g\right>.
\end{equation*}

Zonal (rotation-invariant) functions are those depending only on the first hyperspherical coordinate $\vartheta=\vartheta_1$. Unless it leads to misunderstandings, we identify them with functions of~$\vartheta$ or $t=\cos\vartheta$.
A zonal $\mathcal L^2$-function~$f$ has the following Gegenbauer expansion
\begin{equation}\label{eq: Gegenbauer expansion}
f(t)=\sum_{l=0}^\infty\widehat f(l)\,C_l^\lambda(t),\qquad t=\cos\vartheta,
\end{equation}
where $\widehat f(l)$ are the Gegenbauer coefficients of~$f$ and $\lambda$ is related to the space dimension by
$$
\lambda=\frac{n-1}{2}.
$$
Consequently, for a zonal $\mathcal L^2$-function $f$ one has
$$
\widehat f(l)=A_l^0\cdot a_l^0(f),
$$
compare~\eqref{eq: Fs}, \eqref{eq: Ylk}, and~\eqref{eq: Gegenbauer expansion}.

For $f,g\in\mathcal L^1(\mathcal S^n)$, $g$ zonal, their convolution $f\ast g$ is defined by~\cite[Definition~2.1.1]{DX13}
$$
(f\ast g)(x)=\frac{1}{\Sigma_n}\int_{\mathcal S^n}f(y)\,g(x\cdot y)\,d\sigma(y).
$$
With this notation one has
$$
f\in\mathcal H_l\Longrightarrow f=\frac{\lambda+l}{\lambda}\left(f\ast C_l^\lambda\right)
$$
(Funck-Hecke formula), i.e., the function $\frac{\lambda+l}{\lambda}\, C_l^\lambda$ is the reproducing kernel for~$\mathcal H_l$.

The Poisson kernel on the sphere is given by
\begin{equation}\label{eq: Poisson kernel}\begin{split}
p_r^\lambda(y)&=\frac{1}{\Sigma_n}\sum_{l=0}^\infty r^l\,\frac{\lambda+l}{\lambda}\,\mathcal C_l^\lambda(\cos\vartheta)\\
&=\frac{1}{\Sigma_n}\frac{1-r^2}{(1-2r\cos\vartheta+r^2)^{(n+1)/2}},\quad r=e^{-\rho},\,y_1=\cos\vartheta.
\end{split}\end{equation}
Since the Gegenbauer polynomials~$\mathcal C_l^\lambda$ over the interval $[-1,1]$ are bounded by
\begin{equation*}
\lvert\mathcal C_l^\lambda(\cos\vartheta)\rvert\leq(n+l-2)^{n-2}
\end{equation*}
uniformly in $\vartheta\in[0,\pi]$ (compare \cite[Theorem~7.33.1]{Sz75}), the series in Eq.~\eqref{eq: Poisson kernel} is absolutely convergent for $r\in[0,1)$.

\section{Biharmonic equation}\label{sec: biharmonic equation}

The aim of the present paper is to derive explicit formulas of the Green function~$G_a$ to the biharmonic equation
$$
\left((\Delta_{\mathcal S^n})^2+a\right)u=f,\qquad a\in\mathbb R,\,u\in\mathcal C^2(\mathcal S^n),\,f\in\mathcal C(\mathcal S^n).
$$
Its series representation can be found in many publications, especially for the case $a=0$, see, e.g., \cite[Chapter~5]{FS22}, however, convergence proof is quite often missing. In the proof of the following theorem, considerations concerning convergence are based on the spherical wavelet theory.

\begin{thm}\label{thm:PDEsolution}

\begin{enumerate}

\item
Suppose that $f\in\mathcal C(\mathcal S^n)$ and $u\in\mathcal C^2(\mathcal S^n)$ satisfy
\begin{equation*}
(\Delta_{\mathcal S^n})^2 u+au=f
\end{equation*}
for some $a\in\mathbb R\setminus\{-L^2(L+2\lambda)^2,\,L\in\mathbb N_0\}$. Then,
\begin{equation*}
u=f\ast G_a,
\end{equation*}
for
\begin{equation}\label{eq: Green function series all l}
G_a=\sum_{l=0}^\infty\frac{1}{l^2(l+2\lambda)^2+a}\,\frac{\lambda+l}{\lambda}\,\mathcal C_l^\lambda.
\end{equation}

\item Let $a=-L^2(L+2\lambda)^2$ for some $L\in\mathbb N_0$ and suppose that $f\in\mathcal C(\mathcal S^n)$ and $u\in\mathcal C^2(\mathcal S^n)$ satisfy
\begin{equation*}
(\Delta_{\mathcal S^n})^2 u+au=f,\qquad f\ast\mathcal C_L^\lambda=0,\qquad u\ast\mathcal C_L^\lambda=0.
\end{equation*}
Then,
\begin{equation*}
u=f\ast G_a,
\end{equation*}
for
\begin{equation}\label{eq: Green function series l ne L}
G_a=\sum_{l=0,l\ne L}^\infty\frac{1}{l^2(l+2\lambda)^2+a}\,\frac{\lambda+l}{\lambda}\,\mathcal C_l^\lambda.
\end{equation}

\end{enumerate}

\end{thm}

{\bf Proof.} Analogous to the proof of~\cite[Theorem~1]{INS23}.\hfill\begin{Large}\begin{Large}\smiley\end{Large}\end{Large}

\begin{rem}\label{rem: Ga + const Cllambda} Since in case~(2) it is supposed that $\widehat f(L)=\widehat u(L)=0$, $G_a$ in~\eqref{eq: Green function series l ne L} can be replaced by any function of the form
$$
\sum_{l=0,l\ne L}^\infty\frac{1}{l^2(l+2\lambda)^2+a}\,\frac{\lambda+l}{\lambda}\,\mathcal C_l^\lambda+\text{const}\cdot\mathcal C_L^\lambda.
$$
\end{rem}

\section{Integral representations of the Green function}\label{sec: integral representation}
For $a\in[-1,0]$ the series representations of~$G_a$ in~\eqref{eq: Green function series all l} and~\eqref{eq: Green function series l ne L} can be replaced by an integral representation. It will be shown in the next section how to compute those integrals explicitly.

The first case we consider in $a\in[-1,0)$. It is essentially different from the case $a=0$ and therefore, will be treated separately.

\begin{thm}\label{thm: Green function series Lne0}
\label{thm: Green function series} Let $n\in\mathbb N$, $n\geq2$, be fixed, $\lambda=\frac{n-1}{2}$, and suppose that 
\begin{equation}\label{eq: a vs L}
a=-L^2(L+2\lambda)^2
\end{equation}
for some $L\in(0,\infty)$ such that $\sqrt{\lambda^2-2L\lambda-L^2}\in\mathbb R$, i.e., 
\begin{equation}\label{eq: constraint L}
0<L\leq(\sqrt2-1)\cdot\lambda.
\end{equation}
Denote by~$G_a$ the function
\begin{equation}\label{eq: Green function integer L}
G_a:=\sum_{l=0,l\ne L}^\infty\frac{1}{l^2(n+l-1)^2+a}\frac{\lambda+l}{\lambda}\,\mathcal C_l^\lambda\qquad\text{if }L\in\mathbb N,
\end{equation}
respectively
\begin{equation}\label{eq: Green fuction noninteger L}
G_a:=\sum_{l=0}^\infty\frac{1}{l^2(n+l-1)^2+a}\frac{\lambda+l}{\lambda}\,\mathcal C_l^\lambda\qquad\text{if }L\notin\mathbb N.
\end{equation}
Then
\begin{align}
G_a(\cos\vartheta)&=\int_0^1\left(A\cdot r^{-L-1}-A\cdot r^{L+2\lambda-1}
   +B\cdot r^{\lambda+\sqrt{\lambda^2-2L\lambda-L^2}-1}-B\cdot r^{\lambda-\sqrt{\lambda^2-2L\lambda-L^2}-1}\right)\notag\\
&\cdot\left[\Sigma_n\cdot p_r^\lambda(\cos\vartheta)-\sum_{l=0}^{\lfloor L\rfloor} r^l\cdot\frac{\lambda+l}{\lambda}\,\mathcal C_l^\lambda(\cos\vartheta)\right]\,dr\label{eq: GL}\\
&+\sum_{l=0}^{\lceil L\rceil-1}\frac{1}{l^2(n+l-1)^2+a}\cdot\frac{\lambda+l}{\lambda}\,\mathcal C_l^\lambda(\cos\vartheta)\notag
\end{align}
with
\begin{equation}\label{eq: coefficients A B}
A=\frac{1}{4 L (L + \lambda) (L + 2 \lambda)},\qquad B=\frac{1}{4 L (L + 2 \lambda) \sqrt{\lambda^2-2L\lambda-L^2})}.
\end{equation}
\end{thm}

\begin{rem}\begin{enumerate}
\item For a given $a<0$, $L$ can be chosen to be equal to $\sqrt{\lambda^2+\sqrt{-a}}-\lambda$ such that~\eqref{eq: a vs L} is satisfied.
\item For~$L$ satisfying~\eqref{eq: constraint L}, the number~$a$ is in the interval~$[-1,0)$.
\end{enumerate}
\end{rem}

{\bf Proof. }For $l\ne L$ one has the following decomposition:
\begin{align}
&\frac{1}{l^2(l+2\lambda)^2+a}\notag\\
&\equiv\frac{A}{l-L}-\frac{A}{l+L+2\lambda}+\frac{B}{l+\lambda+\sqrt{\lambda^2-2L\lambda-L^2}}-\frac{B}{l+\lambda-\sqrt{\lambda^2-2L\lambda-L^2}}\label{eq: fraction decomposition}
\end{align}
Further, for $k>-l$,
$$
\int_0^1 r^{l+k-1}\,dr=\frac{1}{l+k}.
$$
Now, according to~\eqref{eq: Poisson kernel}, the expression in brackets on the right-hand-side of~\eqref{eq: GL} equals
\begin{equation}\label{eq: integral r^(l+k-1)}
\sum_{l=\lfloor L\rfloor+1}^\infty r^l\cdot\frac{\lambda+l}{\lambda}\,\mathcal C_l^\lambda(\cos\vartheta)
\end{equation}
and is absolutely convergent for $r\in[0,1)$. Therefore, one can change the order of summation and integration,
\begin{align*}
I_k&:=\int_0^1 r^{k-1}\sum_{l=\lfloor L\rfloor+1}^\infty r^l\cdot\frac{\lambda+l}{\lambda}\,\mathcal C_l^\lambda(\cos\vartheta)\,dr
   =\int_0^1\sum_{l=\lfloor L\rfloor+1}^\infty r^{l+k-1}\cdot\frac{\lambda+l}{\lambda}\,\mathcal C_l^\lambda(\cos\vartheta)\,dr\\
&=\sum_{l=\lfloor L\rfloor+1}^\infty\int_0^1 r^{l+k-1}\cdot\frac{\lambda+l}{\lambda}\,\mathcal C_l^\lambda(\cos\vartheta)\,dr.
\end{align*}
According to~\eqref{eq: integral r^(l+k-1)}, $I_k$ is equal to
$$
I_k=\sum_{l=\lfloor L\rfloor+1}^\infty\frac{1}{l+k}\cdot\frac{\lambda+l}{\lambda}\,\mathcal C_l^\lambda(\cos\vartheta).
$$
This formula can be applied to the integral in~\eqref{eq: GL} with $k=-L$, $k=L+2\lambda$, $k=\lambda+\sqrt{\lambda^2-2L\lambda-L^2}$, and $k=\lambda-\sqrt{\lambda^2-2L\lambda-L^2}$ (note that all these numbers satisfy $k<\lfloor L\rfloor+1$). By fraction decomposition~\eqref{eq: fraction decomposition}, the right-hand-side of~\eqref{eq: GL} equals
\begin{equation}\label{eq: GLseries}
\left(\sum_{l=0}^{l=\lceil L\rceil-1}+\sum_{l=\lfloor L\rfloor+1}^\infty\right)\frac{1}{l^2(l+2\lambda)^2+a}\cdot\frac{\lambda+l}{\lambda}\,\mathcal C_l^\lambda(\cos\vartheta).
\end{equation}
If $L\in\mathbb N$, then $L=\lfloor L\rfloor=\lceil L\rceil$ and the summand with index $l=L$ does not occur in the series~\eqref{eq: GLseries}, and it is equal to the series~\eqref{eq: Green function integer L}. Otherwise, $\lfloor L\rfloor=\lceil L\rceil-1$ and the series~\eqref{eq: GL} is the same as~\eqref{eq: Green fuction noninteger L}.

\hfill\begin{Large}\begin{Large}\begin{Large}\smiley\end{Large}{}\end{Large}\end{Large}

\begin{rem}According to Remark~\ref{rem: Ga + const Cllambda}, $G_a$ can be also expressed as
\begin{align*}
G_a&(\cos\vartheta)= A\cdot\int_0^1 r^{-L-1}
   \cdot\left[\Sigma_n\, p_r^\lambda(\cos\vartheta)-\sum_{l=0}^L r^l\cdot\frac{\lambda+l}{\lambda}\,\mathcal C_l^\lambda(\cos\vartheta)\right]\,dr\\
&+A\cdot\sum_{l=0}^{L-1}\frac{1}{l-L}\cdot\frac{\lambda+l}{\lambda}\,\mathcal C_l^\lambda(\cos\vartheta)\\
&+\int_0^1\left(-A\cdot r^{L+2\lambda-1}+B\cdot r^{\lambda+\sqrt{\lambda^2-2L\lambda-L^2}-1}-B\cdot r^{\lambda-\sqrt{\lambda^2-2L\lambda-L^2}-1}\right)
   \cdot\Sigma_n\, p_r^\lambda(\cos\vartheta)\,dr,
\end{align*}
i.e., splitting of the series representing~$p_r^\lambda$ is performed only for the summand $r^{-L-1}p_r^\lambda$ of the integrand. It is necessary in order to ensure convergence of the integral. For the other exponents: $k=l+2\lambda+1$ and $k=\lambda\pm\sqrt{\lambda^2-2L\lambda-l^2}-1$ the integral $\int_0^1 r^kp_r^\lambda\,dr$ is convergent, however, it is not orthogonal to $\mathcal C_l^\lambda$.
\end{rem}

In the following theorem, an integral representation of the Green function for the case of~$a=0$ is derived.

\begin{thm}\label{thm: G0 integral}
\label{thm: Green function series L=0} Let $n\in\mathbb N$, $n\geq2$, be fixed and let
\begin{equation*}
G_0:=\sum_{l=1}^\infty\frac{1}{l^2(l+2\lambda)^2}\frac{\lambda+l}{\lambda}\,\mathcal C_l^\lambda.
\end{equation*}
Then
\begin{equation}\label{eq: G0 integral}
G_0(\cos\vartheta)=\frac{1}{4\lambda^2}(J_0+J_{2\lambda})-\frac{1}{4\lambda^3}(I_0-I_{2\lambda})
\end{equation}
for
\begin{align*}
I_k&=\int_0^1 r^{k-1}\left[\Sigma_n\cdot p_r^\lambda(\cos\vartheta)-1\right]\,dr,\\
J_k&=\int_0^1\left(\frac{1}{R}\int_0^Rr^{k-1}\left[red\Sigma_n\cdot p_r^\lambda(\cos\vartheta)-1\right]\,dr\right)dR.
\end{align*}
\end{thm}

{\bf Proof. }It was shown in the proof of Theorem~\ref{thm: Green function series Lne0} that
\begin{equation}\label{eq: explicit Ik}
I_k=\sum_{l=1}^\infty\frac{1}{l+k}\cdot\frac{\lambda+l}{\lambda}\,\mathcal C_l^\lambda(\cos\vartheta).
\end{equation}
Further, since for $k>-l$
$$
\frac{1}{(l+k)^2}=\int_0^1\left(\frac{1}{R}\int_0^Rr^{k-1}\cdot r^{l}\,dr\right)dR
$$
holds, we obtain
\begin{equation}\label{eq: explicit Jk}
J_k=\sum_{l=1}^\infty\frac{1}{(l+k)^2}\cdot\frac{\lambda+l}{\lambda}\,\mathcal C_l^\lambda(\cos\vartheta)
\end{equation}
with similar convergence arguments as in the proof of the previous theorem. 

Now, partial fraction decomposition of~$\frac{1}{l^2(l+2\lambda)^2}$ is given by
$$
\frac{1}{l^2(l+2\lambda)^2}=\frac{1}{4\lambda^2l^2} - \frac{1}{4\lambda^3l} + \frac{1}{4\lambda^2(l+2\lambda)^2} + \frac{1}{4\lambda^3(l+2\lambda)}
$$
and this, togehter with~\eqref{eq: explicit Ik} and~\eqref{eq: explicit Jk} yields~\eqref{eq: G0 integral}.\hfill\begin{Large}\begin{Large}\smiley\end{Large}\end{Large}

\section{Explicit expressions of the Green function}\label{sec: explicit representations}

In some cases, for a given set of parameters, the integrals in the representation of~$G_a$ can be computed explicitly. We start the section with some considerations on the exponents $\lambda\pm\sqrt{\lambda^2-2L\lambda-L^2}$, occurring in~\eqref{eq: GL}.

\subsection{Rational exponents $\lambda\pm\sqrt{\lambda^2-2L\lambda-L^2}$}\label{subs: rational exponents}

\begin{lem}\label{lem:rational squares}
Let $\lambda\in\frac12\,\mathbb N$, $L\in\mathbb Q_+\cup\{0\}$. Then,
$$
\sqrt{\lambda^2-2L\lambda-L^2}\in\mathbb Q
$$
if and only if
$$
L=\frac{2(\alpha-1)}{\alpha^2+1}\cdot\lambda\qquad\text{for }\alpha\in\mathbb Q\cap[1,\infty).
$$
\end{lem}

{\bf Proof. }Let $L=\frac{p}{q}\cdot\lambda$, $p\in\mathbb N_0$, $q\in\mathbb N$, gcd$(p,q)=1$. Then
$$
\sqrt{\lambda^2-2L\lambda-L^2}=\frac{\sqrt{2q^2-(p+q)^2}}{q}\cdot\lambda.
$$
This number is real for
\begin{equation}\label{eq: inequality for p L positive}
0\leq p\leq(\sqrt2-1)\cdot q.
\end{equation}
Now, $\sqrt{2q^2-(p+q)^2}\in\mathbb N_0$ if there exists $k\in\mathbb N_0$ such that
$$
2q^2-(p+q)^2=k^2,
$$
i.e.
$$
(q-k)(q+k)=p(p+2q).
$$
Let
$$
\alpha=\frac{q-k}{p}=\frac{p+2q}{q+k}.
$$
From this equation one obtains
$$
q=\frac{(q-k)+(q+k)}{2}=\frac12\left[\alpha p+\frac{p+2q}{\alpha}\right],
$$
and therefore,
\begin{equation}\label{eq: relation p q alpha}
2(\alpha-1)q=(\alpha^2+1)p.
\end{equation}
Equation~\eqref{eq: relation p q alpha} yields
$$
\frac{p}{q}=\frac{2(\alpha-1)}{\alpha^2+1},\qquad\alpha\in\mathbb Q.\\
$$
Further, it can be shown that
$$
0\leq\frac{2(\alpha-1)}{\alpha^2+1}\leq\sqrt2-1
$$
(compare~\eqref{eq: inequality for p L positive}) is satisfied for
$$
\alpha\geq1.
$$
\hfill\begin{Large}\smiley\end{Large}

\subsection{Integer exponents $\lambda\pm\sqrt{\lambda^2-2L\lambda-L^2}$}\label{subs: integer exponents}

Now, note that in the proof of Lemma~\ref{lem:rational squares}, $p$ is an even number. In order to see that, suppose, $\alpha=\frac{a}{b}$, gcm$(a,b)=1$. Then
\begin{equation}\label{eq: p over q in a b}
\frac{p}{q}=\frac{2(\alpha-1)}{\alpha^2+1}=\frac{2b(a-b)}{a^2+b^2}.
\end{equation}
If $a$ and $b$ are of different parity, the statement is obvious. If both $a$ and $b$ are odd, then
$$
a\equiv\pm1\mod4\qquad\text{and}\qquad b\equiv\pm1\mod4.
$$
Therefore
$$
a^2+b^2\equiv2\mod4\qquad\text{and}\qquad 2(a-b)\equiv0\mod4.
$$
Thus, the numerator~$p$ of the quotient in~\eqref{eq: p over q in a b} (when reduced) must be an even number (and the denominator~$q$ an odd number).

Consequently, $\sqrt{2q^2-(p+q)^2}$ is an odd number. Therefore, the numbers $\lambda\pm\sqrt{\lambda^2-2L\lambda-L^2}$ appearing in the fraction denominators on the right-hand-side of~\eqref{eq: fraction decomposition} are equal to
$$
\lambda\pm\sqrt{\lambda^2-2L\lambda-L^2}=\frac12\cdot\left[(n-1)\pm\frac{n-1}{q}\cdot\sqrt{2q^2-(p+q)^2}\right]
$$ 
and they are integer numbers if
\begin{equation}\label{eq: q divisor n-1}
q|(n-1).
\end{equation}
In this case, $L=\frac{p}{q}\cdot\frac{n-1}{2}$ is also an integer. 

Thanks to the constraints~\eqref{eq: q divisor n-1} together with~\eqref{eq: inequality for p L positive}, as well as the fact that \mbox{$2\mid p$}, \mbox{$2\nmid q$}, and $2q^2-(p+q)^2\in\mathbb N^2$, for a given~$\lambda$ one can find all the pairs $(\lambda,L)$ for which $\sqrt{\lambda^2-2L\lambda-L^2}$ is an integer. Table~\ref{tab: pairs lambda L} lists all such pairs with~$L$ different from~$0$ for \mbox{$2\leq n\leq16$}.

\begin{table}\caption{Pairs $(\lambda,L)$, $L\ne0$, that yield integer $\lambda\pm\sqrt{\lambda^2-2L\lambda-L^2}$}\label{tab: pairs lambda L}
\vspace{0.5em}\centering\begin{tabular}{|c|c|c|c|}
\hline
$\lambda$&$L$&$\lambda+\sqrt{\lambda^2-2L\lambda-L^2}$&$\lambda-\sqrt{\lambda^2-2L\lambda-L^2}$\\\hline
$\frac52$ & $1$ & $3$ & $2$\\
$5$ & $2$ & $6$ & $4$\\
$\frac{13}{2}$ & $2$ & $10$ & $3$\\
$\frac{15}{2}$ & $3$ & $9$ & $6$\\\hline
\end{tabular}\end{table}

\begin{rem}There exist pairs of numbers $(\lambda,L)$ with $L<0$ that yield integer values of $\lambda\pm\sqrt{\lambda^2-2L\lambda-L^2}$, different from those listed in Table~\ref{tab: pairs lambda L}. An example is $\lambda=\frac{5}{2}$, $L=-2$ with $\sqrt{\lambda^2-2L\lambda-L^2}=\frac{7}{2}$. In such cases, one obtains a negative~$a$ and there exists a pair $(\lambda,L)$ for which $L(L+2\lambda)=|a|$ that can be found with the procedure described above. 
\end{rem}

\subsection{Green function $G_a$ for $a\in[-1,0)$ on even-dimensional spheres}

If $n$ is even and $L$ as well as $\lambda\pm\sqrt{\lambda^2-2L\lambda-L^2}$ are integer, the integrand in~\eqref{eq: GL} is a rational function of~$r$ and $\sqrt{1-2r\,\cos\vartheta}$. In this case, the function can be integrated with Euler substitution.

In other cases, $G_a$ can be expressed in terms of special functions~\eqref{eq: G in Appell}.

\subsection{Green function $G_a$ for $a\in[-1,0)$ on odd-dimensional spheres}

In the case of an odd dimension~$n$ the Poisson kernel $p_r^\lambda$ is a rational function. If only $L$ and $\lambda\pm\sqrt{\lambda^2-2L\lambda-L^2}$ are rational, the integrals in~\eqref{eq: GL} can be reduced to (or already are) integrals of rational functions. 

\begin{exa} Let $n=3$, i.e., $\lambda=1$, and $L=\frac25$. In this case,
$$
a=-\frac{576}{624}\qquad\text{and}\qquad\sqrt{\lambda^2-2L\lambda-L^2}=\frac15.
$$
The coefficients~$A$ and~$B$ (see formula~\eqref{eq: coefficients A B}) are equal to
$$
A=\frac{125}{672}\qquad\text{and}\qquad B=\frac{125}{96}
$$
and the integrand in~\eqref{eq: GL} is given by
$$
\mathcal I_{1,\frac25}=\frac{125\,(-1+7r^{6/5}-7r^{8/5}+r^{14/5})(-4\cos\vartheta+3r+4r\cos^2\!\vartheta-4r^2\cos\vartheta+r^3)}{672\,r^{2/5}(1-2r\cos\vartheta+r^2)^2}.
$$
Substitute $\rho^5$ for~$r$ to obtain (with computer algebra system Wolfram Mathematica 11.3)
\begin{align*}
\int&\mathcal I_{1,\frac25}\,dr
   =\frac{625}{672}\int\frac{\rho^2(1-7\rho^6+7\rho^8-\rho^{14})(4\cos\vartheta-3\rho^5-4\rho^5\cos^2\vartheta+4\rho^{10}\cos\vartheta-\rho^{15})}
      {(1-2c\rho^5+\rho^{10})^2}\,d\rho\\
&=\frac{625}{672}\left[\frac{\rho^2}{2}+\frac{7\rho^4}{4}-\frac{7\rho^6}{6}+\frac{\rho^{12}}{12}
   +\frac{\rho^2(1-\rho^2)(1+2\rho\cos\vartheta-6\rho^2+2\rho^3\cos\vartheta+\rho^4)}{5(1-2\rho^5\cos\vartheta+\rho^{10})}\right.\\
&\left.-\frac{7}{50}\sum_{\kappa=0}^9\frac{(1-2\rho_\kappa\cos\vartheta+\rho_\kappa^2)(1+\rho_\kappa^4)\ln(\rho-\rho_\kappa)}{\rho_\kappa^3\cos\vartheta-\rho_\kappa^8}\right],
\end{align*}
where~$\rho_\kappa$, $\kappa=0,\dots,9$, are the roots of polynomial $1-2\rho^5\cos\vartheta+\rho^{10}$. Thus,
\begin{align*}
\int_0^1&\mathcal I_{1,\frac25}\,dr=\frac{625}{672}\left[\frac76
   -\frac{7}{50}\sum_{\kappa=0}^9
      \frac{(1-2\rho_\kappa\cos\vartheta+\rho_\kappa^2)(1+\rho_\kappa^4)\ln\left(1-\frac{1}{\rho_\kappa}\right)}{\rho_\kappa^3\,(\cos\vartheta-\rho_\kappa^5)}\right]
\end{align*}
with $\mathcal{I}\left(\ln\left(1-\frac{1}{\rho_\kappa}\right)\right)\in(-\pi,\pi)$. The constant $\frac{625}{672}\cdot\frac{7}{6}$ is equal to~$-\frac{1}{a}$. Consider the terms containing logarithm. Since
\begin{align*}
1-2\rho^5\cos\vartheta+\rho^{10}&=\left(1-\rho^5e^{-i\vartheta}\right)\left(1-\rho^5e^{i\vartheta}\right)=\left(e^{i\vartheta}-\rho^5\right)\left(e^{-i\vartheta}-\rho^5\right)\\
&=\prod_{k=0}^4\left(1-\rho e^{i\frac{\vartheta+2k\pi}{5}}\right)\left(1-\rho e^{-i\frac{\vartheta+2k\pi}{5}}\right),
\end{align*}
the $\rho_\kappa$ are the $5^\text{th}$ roots of $e^{i\vartheta}$ and of~$e^{-i\vartheta}$,
\begin{align*}
\rho_\kappa&=e^{i\frac{\vartheta+2\kappa\pi}{5}},\qquad\,\,\,\kappa=0,\dots,4,\\
\rho_\kappa&=e^{-i\frac{\vartheta+2\kappa\pi}{5}},\qquad\kappa=5,\dots,9.
\end{align*}
For $\kappa=0,\dots,4$ set
$$
\alpha=\frac{\vartheta+2\kappa\pi}{5}
$$
and
$$
v_\kappa:=\frac{(1-2\rho_\kappa\cos\vartheta+\rho_\kappa^2)(1+\rho_\kappa^4)\ln\left(1-\frac{1}{\rho_\kappa}\right)}{\rho_\kappa^3\,(\cos\vartheta-\rho_\kappa^5)}.
$$
Then, using relation $\cos\vartheta=\cos(5\alpha)=\cos\left(5(-\alpha)\right)$, one obtains with Mathematica
\begin{align*}
v_\kappa&+v_{\kappa+5}
   =\frac{4i\left[\ln(-e^{-i\alpha})-\ln(-e^{i\alpha})+\ln(1-e^{i\alpha})-\ln(1-e^{-i\alpha})\right](\sin2\alpha+\sin4\alpha+\sin6\alpha)}{1+2\cos2\alpha+2\cos4\alpha}.
\end{align*}
Further,
$$
\ln(-e^{-i\alpha})-\ln(-e^{i\alpha})+\ln(1-e^{i\alpha})-\ln(1-e^{-i\alpha})=\ln(-e^{-i\alpha})=i(\pi-\alpha)
$$
and, consequently,
$$
v_\kappa+v_{\kappa+5}=\frac{4\left(\alpha-\pi\right)(\sin2\alpha+\sin4\alpha+\sin6\alpha)}{1+2\cos2\alpha+2\cos4\alpha}.
$$
Therefore,
$$
G_a(\cos\vartheta)=\frac{25}{48}\sum_{\alpha\in\mathcal A}\frac{(\pi-\alpha)(\sin2\alpha+\sin4\alpha+\sin6\alpha)}{1+2\cos2\alpha+2\cos4\alpha}
$$
with $\mathcal A=\left\{\frac{\vartheta+2\kappa\pi}{5}:\,\kappa=0,\dots,4\right\}$.
\end{exa}

Otherwise, suppose, $L$ or $\sqrt{\lambda^2-2L\lambda-L^2}$ is irrational. The integrand in~\eqref{eq: GL} is a sum of (a polynomial and) expressions of the form
$$
\frac{r^\alpha}{(1-2r\cos\vartheta+r^2)^{\lambda+1}},\qquad\alpha>0.
$$
According to~\cite[9.3(4)]{wB64},
\begin{equation}\label{eq: Appell_as_integral}\begin{split}
\int_0^1&u^{\alpha-1}(1-u)^{\gamma-\alpha-1}(1-ux)^{-\beta}(1-uy)^{-\beta'}du\\
&=\frac{\Gamma(\gamma-\alpha)\Gamma(\alpha)}{\Gamma(\gamma)}F_1(\alpha;\beta,\beta';\gamma;x,y),
\end{split}\end{equation}
where~$F_1$ is the Appell $F_1$-function and $0<\text{Re }\alpha<\text{Re }\gamma$. Since
$$
1-2r\cos\vartheta+r^2=(1-re^{i\vartheta})(1-re^{-i\vartheta}),
$$
we obtain from~\eqref{eq: Appell_as_integral}
\begin{equation}\label{eq: G in Appell}
\int_0^1\frac{r^\alpha}{(1-2r\cos\vartheta+r^2)^{\lambda+1}}\,dr=\frac{F_1(\alpha;\lambda+1,\lambda+1;\alpha+1;e^{i\vartheta},e^{-i\vartheta})}{\alpha}.
\end{equation}

\subsection{Green function $G_0$ on even-dimensional spheres}

In the case of even-dimensional spheres, I was not able to find a closed formula for~$J_0$ and~$J_{2\lambda}$ defined in~Theorem~\ref{thm: G0 integral}. In Chapter~5 of~\cite{FS22} a closed formula for~$G_0$ on the two-dimensional sphere is derived. The authors solve the differential equation
$$
\Delta_{\mathcal S^n} G(\cos\vartheta)=K(\cos\vartheta),
$$
with properly chosen boundary conditions. $K$ denotes the Green function to the Poisson equation.

\subsection{Green function $G_0$ on odd-dimensional spheres}

\begin{exa}
Let $n=3$. Then, $\lambda=1$ and
$$
f(r):=\Sigma_n\cdot p_r^\lambda(\cos\vartheta)-1=\frac{4r\cos\vartheta-3r^2-4r^2\cos^2\vartheta+4r^3\cos\vartheta-r^4}{(1-2r\cos\vartheta+r^2)^2}.
$$
Thus (with Wolfram Mathematica 11.3),
$$
\int\frac{f(r)}{r}\,dr=\frac{1}{1-2r\cos\vartheta+r^2}+\cot\vartheta\cdot\arctan\frac{r-\cos\vartheta}{\sin\vartheta}-\frac12\ln(1-2r\cos\vartheta+r^2)+C
$$
and
\begin{align*}
I_0&=\int_0^1\frac{f(r)}{r}\,dr=-1+\frac{1}{2-2\cos\vartheta}+\cot\vartheta\cdot\arctan\sqrt\frac{1+\cos\vartheta}{1-\cos\vartheta}-\frac12\ln(2-2\cos\vartheta)\\
&=-1+\frac{1}{2-2\cos\vartheta}+\frac{\pi-\vartheta}{2}\cdot\cot\vartheta-\frac12\ln(2-2\cos\vartheta).
\end{align*}
Similarly,
$$
I_2=\frac{\cos\vartheta}{2-2\cos\vartheta}-\frac{\pi-\vartheta}{2}\cdot\cot\vartheta-\frac12\ln(2-2\cos\vartheta).
$$
Further,
\begin{align*}
J_0&=\int_0^1\frac{1}{R}\left[\frac{2R\cos\vartheta-R^2}{1-2R\cos\vartheta+R^2}+\cot\vartheta\cdot\arctan\left(\frac{R\sin\vartheta}{1-R\cos\vartheta}\right)\right.\\
&\left.\qquad-\frac12\ln(1-2R\cos\vartheta+R^2)\right]dR
\end{align*}
and
\begin{align*}
J_2&=\int_0^1\frac{1}{R}\left[1-\frac{R^2}{2}-\frac{1-2R\cos\vartheta}{1-2R\cos\vartheta+R^2}-\cot\vartheta\cdot\arctan\left(\frac{R\sin\vartheta}{1-R\cos\vartheta}\right)\right.\\
&\left.\qquad-\frac12\ln(1-2R\cos\vartheta+R^2)\right]dR.
\end{align*}
Thus,
\begin{align*}
J_0+J_2&=\int_0^1\left[\frac{4\cos\vartheta-R+2R^2\cos\vartheta-R^3}{2(1-2R\cos\vartheta+R^2)}-\frac1R\ln(1-2R\cos\vartheta+R^2)\right]dR\\
&=\left[-\frac{R^2}{4}+2\cot\vartheta\cdot\arctan\frac{R-\cos\vartheta}{\sin\vartheta}+\mathcal L_2(Re^{i\vartheta})+\mathcal L_2(Re^{-i\vartheta})\right]_0^1\\
&=-\frac{1}{4}+(\pi-\vartheta)\cdot\cot\vartheta+\mathcal L_2(e^{i\vartheta})+\mathcal L_2(e^{-i\vartheta}),
\end{align*}
where~$\mathcal L_2$ denotes the dilogarithm function. Now, by~\cite[formula~(3.2)]{lM03}
$$
\mathcal L_2(e^{i\vartheta})+\mathcal L_2(e^{-i\vartheta})=-\frac16\pi^2-\frac12\ln^2(-e^{i\vartheta})=\frac{\pi^2}{3}-\pi\vartheta+\frac{\vartheta^2}{2},
$$
and alltogether
\begin{equation}\label{eq: G0 n3}
G_0(\vartheta)=\frac{3+4\pi^2-12\pi\vartheta+6\vartheta^2}{48}.
\end{equation}
\end{exa}

This result can obtained also via solution of the equation
\begin{equation}\label{eq: diffeq G}
\Delta_{\mathcal S^n}G_0(\vartheta)=K(\vartheta),
\end{equation}
where~$K$ denotes the Green kernel to the Poisson equation~\cite{INS23}. Contrary to the two-dimensional case, solved in~\cite[Subsection~5.7]{FS22}, for an odd~$n$ formula~\eqref{eq: recursion LB theta}, i.e., the one using variable~$\vartheta$ instead of~$t=\cos\vartheta$, should be applied.

\begin{exa}
For $n=3$ the Green function to the Poisson equation~$\Delta_{\mathcal S^n} u=f$ is given by~\cite[Table~1]{INS23}
$$
K(\vartheta)=\frac14+\frac{\vartheta-\pi}{2}\cdot\cot\vartheta.
$$
Since both~$G_0$ and~$K$ are rotation invariant functions, equation~\eqref{eq: diffeq G} reduces to
$$
\frac{1}{\sin^2\vartheta}\frac{\partial}{\partial\vartheta}\left(\sin^2\vartheta\cdot\frac{\partial G_0(\vartheta)}{\partial\vartheta}\right)=\frac14+\frac{\vartheta-\pi}{2}\cdot\cot\vartheta.
$$
Thus,
\begin{align*}
\sin^2\vartheta&\cdot\frac{\partial G_0(\vartheta)}{\partial\vartheta}=\int\left(\frac14\sin^2\vartheta+\frac{\vartheta-\pi}{4}\cdot\sin2\vartheta\right)d\vartheta\\
&=\frac{\vartheta-\pi}{4}\cdot\sin^2\vartheta+C_1,
\end{align*}
and further
$$
G_0(\vartheta)=\int\left(\frac{\vartheta-\pi}{4}+\frac{C_1}{\sin^2\vartheta}\right)d\vartheta=\frac{\vartheta^2}{8}-\frac{\pi\vartheta}{4}+C_1\cot\vartheta+C_2.
$$
Now, since
$$
G_0(0)=\sum_{l=0}^\infty\frac{1}{l^2(l+2)^2}\cdot\frac{l+1}{1}\cdot(l+1)=\frac{3+4\pi^2}{48}<\infty,
$$
the constant~$C_1$ vanishes and $C_2=\frac{3+4\pi^2}{48}$, i.e., $G_0$ is given by~\eqref{eq: G0 n3}.
\end{exa}

\end{document}